\documentclass{amsart}

\begin{document}

\title{A note on inhomogeneous foliations with sections}
\author{Stephan Wiesendorf}

\address{Stephan Wiesendorf\\
Mathematisches Institut, Universit\"at zu K\"oln, Weyertal 86-90, 50931 K\"oln, Germany}
\email{swiesend@math.uni-koeln.de}


\begin{abstract} We give an easy example showing that sections of a singular Riemannian foliation on a simply connected space neither have to be isometric nor injectively immersed.
\end{abstract}
  
\maketitle 

\noindent If a group $G$ acts on a Riemannian manifold via isometries, the partition $\mathcal F^G$ of $M$ into the orbits defines a so called \emph{singular Riemannian foliation} on $M$, i.e. a smooth partition into injectively immersed submanifolds, called the leaves, with the property that a geodesic is orthogonally to the leaves at either all or none of its points. The foliation is said to have \emph{sections} if through every point lying on a leaf of maximal dimension, there is a complete immersed submanifold which meets every orbit and always orthogonally. An isometric action of a group is called \emph{polar} if the corresponding foliation $\mathcal F^G$ admits sections. We refer the reader who is not familiar with these concepts to \cite{GZ} and \cite{A}.

If a foliation admits sections, all sections are locally isometric and have the same universal cover (cf. \cite{T}). For a generic section $\Sigma$, T\"oben constructed in \cite{T} a generalized Weyl group $W(\Sigma)$ acting on $\Sigma$ with $M/\mathcal F = \Sigma/W(\Sigma)$. But it is neither known if there exists such a group for every section nor if there are general (topological) assumptions such that all sections have to be equal. If a group acts polarly on a Riemannian manifold, the group also acts transitively on the set of sections. In particular, all the sections are isometric. That this is generally not true in the inhomogeneous case can easily be observed by the example of the Klein bottle as the nonorientable $S^1$-bundle over $S^1$ viewed as the soul of the M\"obius band identification. In this case, this soul represents an exeptional section which is covered twice by the other sections.

With this example in mind one could still believe that these phenomena are caused by the fact that $M$ is not simply connected. This belief might be motivated by the fact that singular Riemannian foliations on simply connected spaces cannot be as complicated as in the general setting. For instance, in \cite{L} it is shown that if the foliation admits sections, all leaves have to be closed and are therefore embedded in this case. Moreover, due to \cite{AT}, a foliation with sections on a simply connected space cannot have exeptional leaves, i.e. the restriction of the foliation to the open set of all points lying on leaves of maximal dimension is induced by a Riemannian submersion (see \cite{L} for more properties of foliations on simply connected spaces). So, it seems natural to ask if all sections of a singular Riemannian foliation have to be isometric if the manifold is simply connected.  In the following, we give an easy example which shows that, even in the simplest cases, this does not have to be true. In fact, the behavior of the sections of an inhomogeneous foliation can be arbitrary.

Our example also shows that sections are not necessarily injectively immersed, even in the simply connected case. Whereas, the corresponding question in the homogeneous case, stated by Grove and Ziller in \cite{GZ}, is still open, except for the codimension one case where it is known by a result in \cite{AA} that sections are either closed or injectively immersed.   

\vspace{0.5cm}



\noindent \textbf{Example:}
Let $D^n$ denote the closed unit $n$-ball in $\mathbb R^n$ foliated by concentric spheres, i.e. by the one-parameter family $L_t = S^{n-1}(t)$, where $t$ is just the radius. Now, we glue two copies of $D^n$ along their common boundary $L_1= S^{n-1}$ by a diffeomorphism $f : S^{n-1} \to S^{n-1}$ to obtain a new manifold $S = D^n \sqcup_fD^n$, homeomorphic to $S^n$, with a codimension one foliation $\mathcal F$. 

We now choose a metric $g$ on $S$ such that in each copy of $D^n$ the horizontal geodesics are radial segments. This can be done as follows. Around the singularities let $g$ coincide with the standard Euclidean metric. Over $D^n \setminus \{0\}$, the tangent bundle $T\mathbb R^n$ splits  as a direct sum $V \oplus H$, where the distribution $V$ is pointwise given by the tangent spaces of the leaves of $\mathcal F$, and $H$ is just the line bundle $V^{\perp}$ (with respect to the standard metric). Since the gluing respects $V$, there is an obvious gluing of the two copies of $T\mathbb R^n|_{D^n}$ that induces a splitting of $TS$ outside the singularities into $\bar V \oplus \bar H$ such that $\bar V = T\mathcal F$.  By construction, the standard metric induces a metric on $\bar H$ and we declare $\bar V$ and $\bar H$ to be $g$-orthogonal and define $g$ to coincide with the standard metric on $\bar H$. Around the gluing leaf, we choose an arbitrary Riemannian metric on $\bar V$ and interpolate smoothly between this and the standard metric while moving towards the singularities by means of the radius. This gives us a Riemannian metric on $S$ that makes $\mathcal F$ into a singular Riemannian foliation with the desired properties. Namely, the sections are the horizontal geodesics which are piecewise given by radial lines $t \mapsto tp$ for $p \in S^{n-1}$. Thus, starting in a singular point a (unit speed) geodesic $c(t)$ is a straight line for $t \in [0,1]$. Then, one has to identify the point $c(1)$ with $f(c(1))$ and move along the line $<f(c(1))>$ up to the point $-f(c(1))$ which is then identified with $f^{-1}(-f(c(1))$ and one moves along the corresponding line up to $-f^{-1}(-f(c(1))$, and so on. 

Therefore, we see that whether such a geodesic is closed only depends on the diffeomorphism $f : S^{n-1} \to S^{n-1}$. Namely, it is closed iff $(-f^{-1} \circ (-f))^k(c(1)) = c(1)$ for some $k \in \mathbb N$. We call the least number $k$ with this property the $\emph{period}$ of $c(1)$ and say that it has $\emph{no period}$ if there is no such $k$. 

Thus, all we have to show is that there are diffeomorphisms $S^{n-1} \to S^{n-1}$ such that the period is not constant. But those maps are easy to construct in any dimension. Note that a point has period $1$ if and only if $f(-p) = -f(p)$, i.e. $f$ maps the line through $p$ onto the line through $f(p)$. We therefore only have to consider a diffeomorphism $f : S^1 \to S^1$ with the property that $f$ is the identity on one (closed) semicircle and differs from the identity on the other semicircle.

\vspace{0.5cm} 

\noindent \emph{Remarks.} In any dimension, one can choose a diffeomorphism $f$ with the above properties, namely the occurrence of different as well as no periods, such that $S$ is actually diffeomorphic to $S^n$.

In \cite{M} it is incorrectly stated that for a closed singular Riemannian foliation of codimension 1 on a closed manifold all horizontal geodesics have to be closed (cf. Lemma 3.2), what is even not true in the homogeneous case. For example, one could take the diagonal $\mathrm{SO}(3)$-action on the product $S^2(r) \times S^2(s)$ with $r/s$ irrational. Our example mainly describes the phenomena overlooked in the proof given there.

\vspace{0.5cm}

\noindent \textbf{Acknowledgements.} I would like to thank Alexander Lytchak for helpful discussions on that topic and on a previous version of this note and Marcos Alexandrino who brought the paper \cite{M} to my attention. I am also grateful to Wolfgang Ziller for the example at the end of this note.

\end{document}